\documentclass[12pt,dvips]{amsart}
\usepackage{euler, amsfonts, amssymb, latexsym, epsfig,epic, stmaryrd}

\newtheorem{Theorem}{Theorem}[section]
\newtheorem{Proposition}[Theorem]{Proposition} 
\newtheorem{Lemma}[Theorem]{Lemma}

\newtheorem{Corollary}[Theorem]{Corollary}

\newtheorem*{Corollary*}{Corollary}
\newtheorem*{Theorem*}{Theorem}

\newtheorem{thm}{Theorem}

\theoremstyle{remark}
 
\newtheorem{Remark}[Theorem]{Remark}
\newtheorem{Example}[Theorem]{Example}
\newtheorem*{Example*}{Example}

\theoremstyle{plain}

\newcommand\dfn[1]{\textbf{#1}} 

\newcommand\integers{{\mathbb Z}}
\newcommand\naturals{{\mathbb N}}
\newcommand\groth{{\mathfrak G}}

\newcommand\kk{{\Bbbk}}
                  
\newcommand\bxty{{(x \notmapsto y)}}
\newcommand\codim{{\rm codim}}
\newcommand\del{{\rm del}}
\newcommand\link{{\rm link}}
\newcommand\starr{{\rm star}}

\def\<{\langle}
\def\>{\rangle}
\def\fillrightmap{\mathord- \mkern-6mu
        \cleaders\hbox{$\mkern-2mu \mathord- \mkern-2mu$}\hfill
        \mkern-6mu \mathord\rightarrow}
\def\TO#1{\begin{array}{@{}c@{}}\\[-3.5ex]
        \ \scriptstyle #1\;\,\\[-1.4ex]\fillrightmap\end{array}}
\def\DELTA#1#2#3#4{\Delta_{#4}(#1 \,\TO{#3}\, #2)}
\def\K{$K$}
\def\BT{{\bigcup T}}
\def\ID{{I_{\hspace{-.15ex}\Delta}}}
\def\IDstar{{I_{\hspace{-.2ex}\Delta}^\star}}

\def\To{\Rightarrow}
\def\notmapsto{\arrownot\mapsto}
\def\nothing{\varnothing}
\def\mym{{(m \notmapsto y_m)}}
\def\SVT{\mathit{SVT}}
\def\SSYT{\mathit{SSYT}}
\def\LSVT{\mathit{LSVT}}

\newcommand{\comment}[1]{$\star${\sf\textbf{#1}}$\star$}
\newcommand\excise[1]{}

\renewcommand\tt{{\mathbf{t}}}

\newcommand{\cellsize}{12}
\newlength{\cellsz} 
\newsavebox{\cell}
\sbox{\cell}{\begin{picture}(\cellsize,\cellsize)
\put(0,0){\line(1,0){\cellsize}}
\put(0,0){\line(0,1){\cellsize}}
\put(\cellsize,0){\line(0,1){\cellsize}}
\put(0,\cellsize){\line(1,0){\cellsize}}
\end{picture}}
\newcommand\cellify[1]{\def\thearg{#1}\def\nothing{}%
\ifx\thearg\nothing
\vrule width0pt height\cellsz depth0pt\else
\hbox to 0pt{\usebox{\cell} \hss}\fi%
\vbox to \cellsz{
\vss
\hbox to \cellsz{\hss$#1$\hss}
\vss}}
\newcommand\tableau[1]{\vtop{\let\\\cr
\baselineskip -16000pt \lineskiplimit 16000pt \lineskip 0pt
\ialign{&\cellify{##}\cr#1\crcr}}}

\begin{document}
\pagestyle{plain}
\mbox{}
\vspace{2ex}

\title{Tableau complexes}
\author{Allen Knutson} 
\address{Department of Mathematics, University of California,
  San Diego, CA 92093, USA}
\email{allenk@math.ucsd.edu}

\author{Ezra Miller} 
\address{Department of Mathematics, University of Minnesota,
  Minneapolis, MN 55455, USA} 
\email{ezra@math.umn.edu}

\author{Alexander Yong} 
\address{Department of Mathematics, University of Minnesota,
  Minneapolis, MN 55455, USA; and
  \newline\indent The Fields Institute, 222 College Street, Toronto,
  Ontario, M5T 3J1, Canada}
\email{ayong@math.umn.edu} 
\date{17 March 2006}

\begin{abstract}
  Let $X,Y$ be finite sets and $T$ a set of functions from $X\to Y$ which
  we will call ``tableaux''. We define a simplicial complex
  whose facets, all of the same dimension, correspond to these 
  tableaux.  Such {\em tableau complexes} have many nice properties,
  and are frequently homeomorphic to balls, which we prove using
  vertex decompositions~\cite{BP}.

  In our motivating example, the facets are labeled by semistandard
  Young tableaux, and the more general interior faces are labeled by
  Buch's set-valued semistandard tableaux.  One vertex decomposition
  of this ``Young tableau complex'' parallels Lascoux's transition
  formula for vexillary double Grothendieck polynomials \cite{lascoux,
  lascoux:ASM}.  Consequently, we obtain formulae (both old and new)
  for these polynomials.  In particular, we present a common
  generalization of the formulae of Wachs \cite{Wachs} and Buch
  \cite{buch:KLR}, each of which implies the classical tableau formula
  for Schur polynomials.
%
\end{abstract}

\maketitle

\tableofcontents

\section{Introduction}

\subsection{Statement of results}

Let $X$ and~$Y$ be two finite sets. We will call functions from $X$ to~$Y$
~\dfn{tableaux}; we think of each tableau $f:X\to Y$ as a labeling
of the points of~$X$ by elements of~$Y$.
Formally, we identify a tableau $f$ with its corresponding set
$\{(x \mapsto y) : f(x) = y\} \subseteq X \times Y$ of ordered pairs,
whose projection to $X$ is bijective.  

We specify a subset $T$ of ``special'' tableaux. Also, 
let $E \subseteq X \times Y$ be
a relation containing every $f \in T$.  There are obvious minimal and
maximal choices of $E$, namely $\BT := \bigcup_{f \in T} f$ and $X
\times Y$, but it will be convenient to not~restrict~$E$.

Our motivating example is when $X$ is the set of boxes in
partition~$\lambda$ and $Y = \{1,\ldots,n\}$, so a tableau $f : X \to
Y$ is a Young tableau of shape~$\lambda$ and with entries bounded above by
$n$ (without any other demands on the labeling), 
and $T$ is the special subset of semistandard Young tableaux.  In a moment
(Section~\ref{sub:youngtab}), we will describe this case in detail.

Define the simplicial complex $\DELTA XYTE$, which we call a
\dfn{tableau complex}, as follows.  Consider the collection of subsets
of~$E$ as a simplex under reverse inclusion; thus the vertices are the
complements $\bxty := E \setminus \{(x \mapsto y)\}$ rather than the
elements $(x \mapsto y) \in E$.  We view the faces of this simplex as
\dfn{set-valued tableaux}, thought of as relations $F: X \To Y$, in
which every element $x \in X$ is labeled by a set of elements $F(x)
\subseteq Y$.  A set-valued tableau $F'$ is a face of~$F$ whenever $F'
\supseteq F$, meaning that $F'(x) \supseteq F(x)$ for all~$x \in X$.
(This set-theoretic containment is always intended when we say that
one set-valued tableau contains another, even when both are being considered as
faces of a simplicial complex, where containment among faces goes the
opposite way.)
The tableau complex is defined by its facets (maximal faces), which we
declare to be the tableaux $f \in T$.  (In this paper, the terms
``function'' and ``tableau'', when unadorned by ``set-valued'', mean
single-valued functions in the usual sense.)  The face with no
vertices, which we call the \dfn{empty face}, is the set-valued
tableau~$E$.

\begin{Example*}
Consider the tableau complex in which $X=\{1,2,3,4\}$ and $Y$ is the
English alphabet, and where $T$ consists of all the English words in
[dh]e[al][dl].  In detail, $E = \{(1\mapsto d),(1\mapsto h),(2\mapsto
e),(3\mapsto a), (3\mapsto l),(4\mapsto d),(4\mapsto l)\}$.  For
simplicity, in the following figure, at the vertices we indicate the
{\em unused} letters of~$E$.
\begin{figure}[htbp]
  \centering
  \epsfig{file=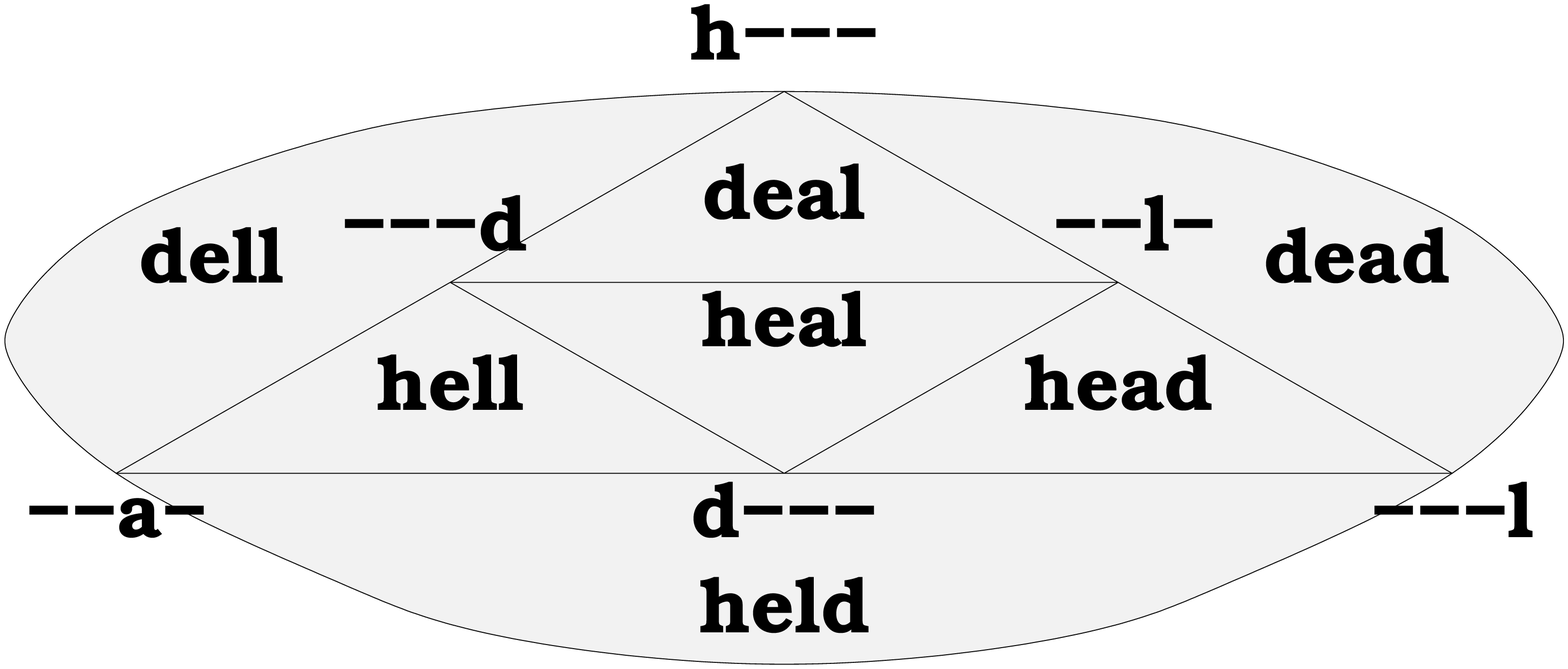,height=1.5in}
  \label{fig:deld}
\end{figure}
For example, the vertex common to ``hell'', ``heal'', ``head'' and
``held'' is $(1 \notmapsto d) = E\setminus \{(1 \mapsto d)\}$.  This
complex is a $2$-ball, but if ``deld'' were a word, it would label the
outer face, and this complex would be a $2$-sphere.
\end{Example*}

\begin{thm}\label{t:tabcomplex}
  The following hold for an arbitrary tableau complex $\Delta = \DELTA
  XYTE$.
  \begin{enumerate}
  \item
    $\Delta$ is \dfn{pure}, meaning that its facets (the tableaux $f
    \in T$) all have the same dimension.
  \item
    The codimension of a face~$F$ in~$\Delta$ is the number $\sum_{x
    \in X} (|F(x)|-1)$ of its ``extra'' values.  For example, faces of
    codimension~$1$, which we call \dfn{ridges}, are set-valued
    tableaux taking two values for precisely one $x \in X$ and taking
    one value at all other points of~$X$.
  \item
    Each ridge is a face of at most two facets.  In particular, if a
    tableau complex is shellable then it is homeomorphic to a ball or
    sphere.
  \item
    The link of a face in a tableau complex is again a tableau
    complex.
  \end{enumerate}
\end{thm}
All of these results are proved in Section~\ref{sub:gen} except the
last (the statement about links), which is
Proposition~\ref{prop:link}.

We shall see that abstractly, a tableau complex is a (multicone on a)
top-dimensional subcomplex of a join of boundaries of simplices.
Tableau complexes can also be characterized, among pure complexes, by
the extremal property given toward the end of
Section~\ref{sec:extremal}.

Although tableau complexes are not generally (shellable) balls or
spheres, we can give conditions that guarantee this conclusion.  The
next theorem thus defines the main class of tableau complexes of
interest in this paper.  Except for the claim about interior faces,
which is Proposition~\ref{prop:boundary}, it is a simpler-to-state
special case of Theorem~\ref{thm:posetball}.

\begin{thm}\label{t:posetball}
  Let $X$ be a poset, and $Y$ totally ordered.
  Let $\Psi$ be a set of pairs $(x_1,x_2)$ in $X$ with $x_1 < x_2$.
  Let $T$ be the set of tableaux $f: X \to Y$ such that
  \begin{itemize}
  \item if $x_1\leq x_2$, then $f(x_1) \leq f(x_2)$, and
  \item if $(x_1,x_2)\in \Psi$, then $f(x_1) < f(x_2)$;
  \end{itemize}
  thus $T$ consists of the order-preserving tableaux from $X$ to~$Y$
  that are strictly order-preserving on the pairs in~$\Psi$.  Let $E
  \supseteq \BT$.  Then the tableau complex $\DELTA XYTE$ is
  \begin{enumerate}
  \item
  homeomorphic to a ball or sphere;
  \item
  vertex-decomposable, as defined in \cite{BP}, and hence shellable; and
  \item
  a manifold with (possibly empty) boundary whose interior faces are
  those set-valued tab\-leaux~$F$ such that every 
  tableau $f \subseteq F$ lies in~$T$.
  \end{enumerate}
\end{thm}

If $Y$ is taken to be a set of natural numbers, then the tableaux in
Theorem~\ref{t:posetball} are {\em $P$-partitions} \cite{EC1}, where
$X=P$.  In our context, however, this point of view is misleading for
a couple of reasons.  First, the condition that $Y$ be totally ordered
can be relaxed in a natural way, as we will see during the proof.
Second, $P$-partitions naturally form a set that is infinite and
possesses additive structure; both of these properties are unnatural
from the point of view of tableau complexes.  More deeply,
$P$-partitions correspond naturally to the basis elements of the
Stanley-Reisner ring of a certain simplicial complex (a Gr\"obner
degeneration of the cone of $P$-partitions) rather than to the facets.

We give three formulae for the Hilbert series of the Stanley-Reisner
ring, the third one based on an explicit shelling of tableau
complexes.  For proofs, see Section~\ref{sec:hilb}, where the
statements break the products over $v \supseteq F$ and $v
\not\supseteq F$ further into products over~$x \in X$.

\begin{thm}\label{t:kpoly}
Let $\Delta = \DELTA XYTE$ be a tableau complex, and recall that the
vertices of~$\Delta$ are set-valued tableaux $\bxty \subseteq E$.
\begin{enumerate}
\item\label{C:1}
The Hilbert series, in variables $\{t_v : v$ is a vertex
of~$\Delta\}$, equals $K_\Delta/\prod (1-t_v)$, where the
denominator product is over all vertices $v$ of~$\Delta$, and the
numerator is the \K-polynomial
\[
  K_\Delta = \sum_F \prod_{v \supseteq F}t_v \prod_{v \not\supseteq
  F}(1-t_v),
\]
the sum being over all set-valued tableaux $F \subseteq E$ such that
$f \subseteq F$ for some $f \in T$.

\item\label{C:2}
If $\Delta$ is homeomorphic to a ball or a sphere, then writing $|F| =
\sum_{x \in X} |F(x)|$ and $|X|$ for the size of~$X$, the
\K-polynomial can be expressed an alternating sum
\[
  K_\Delta = \sum_F (-1)^{|F|-|X|}\prod_{v \not\supseteq F}(1-t_v)
\]
over the set-valued tableaux $F \subseteq E$ such that every tableau
$f \subseteq F$ satisfies $f \in T$.

\item\label{C:3}
Assume furthermore the hypotheses of Theorem~\ref{t:posetball}, and
set $E = \BT$.  Then there is a shelling for~$\Delta$ such that the
minimal new face when the facet $f \in T$ is added during the shelling
is an explicitly described set-valued tableau $N(f) \supseteq f$.
Consequently,
\[
  K_\Delta = \sum_{f\in T} 
  \prod_{v \not\supseteq f} (1-t_{v})
  \prod_{v \supseteq N(f)} t_{v}.
\]
\end{enumerate}
\end{thm}

As a result of Theorem~\ref{t:kpoly}.\ref{C:3}, we get a positive
combinatorial rule to compute the $h$-vector $(h_0, h_1,\ldots \ )$ of
$\Delta$: if $\eta(f) = |E\setminus N(f)|-1$, then $h_j$ counts the
number of $f \in T$ with $\eta(f)=j$.

\subsection{Young tableau complexes}\label{sub:youngtab}

We now describe the prototypical example of a tableau complex, and an
application to computing vexillary Grothendieck polynomials, or
equivalently, Hilbert series formulae for vexillary determinantal
varieties.

Let $\lambda \subseteq \naturals^2$ be an English partition, or
equivalently, a Young shape with its origin at its upper-left corner.
A \dfn{set-valued Young tableau} \cite{buch:KLR} is a filling of the
boxes of~$\lambda$, each with a nonempty finite set of natural
numbers.  The set in each box is typically expressed as a strictly
increasing list.  When the set in every box is a singleton, what
results is an (ordinary) Young tableaux.  If $|\tau|$ denotes the
number of entries in a set-valued tableau~$\tau$, and $|\lambda|$ is
the number of boxes in the partition, then $|\tau| \geq |\lambda|$.
Moreover, $|\tau|=|\lambda|$ only for tableaux.  (Tableaux are assumed
ordinary unless the term ``set-valued'' is written.)

A set-valued tableau $\tau$ is called \dfn{semistandard} if for every pair
$b_1,b_2$ of boxes of $\tau$,
\begin{itemize}
\item
  each entry of $b_1$ is weakly less than each entry of $b_2$ whenever
  $b_1$ lies left of~$b_2$, and
\item
  each entry of $b_1$ is strictly less than each entry of $b_2$
  whenever $b_1$ lies above $b_2$.
\end{itemize}
One can speak of one set-valued tableau \dfn{containing} another (of
the same shape $\lambda$) if for each box of~$\lambda$, the set of
numbers in one set-valued tableau contains the corresponding set in
the other.  In these terms, for $\tau$ to be semistandard, one needs
that every tableau contained in~$\tau$ is semistandard in the usual
sense.  More generally, we define a set-valued tableau to be
\dfn{limit semistandard} if {\em some} tableau it contains is
semistandard.  For example, the first of the following set-valued
tableaux is semistandard, the second is limit semistandard, and the
third is neither:
\[
  \epsfig{file=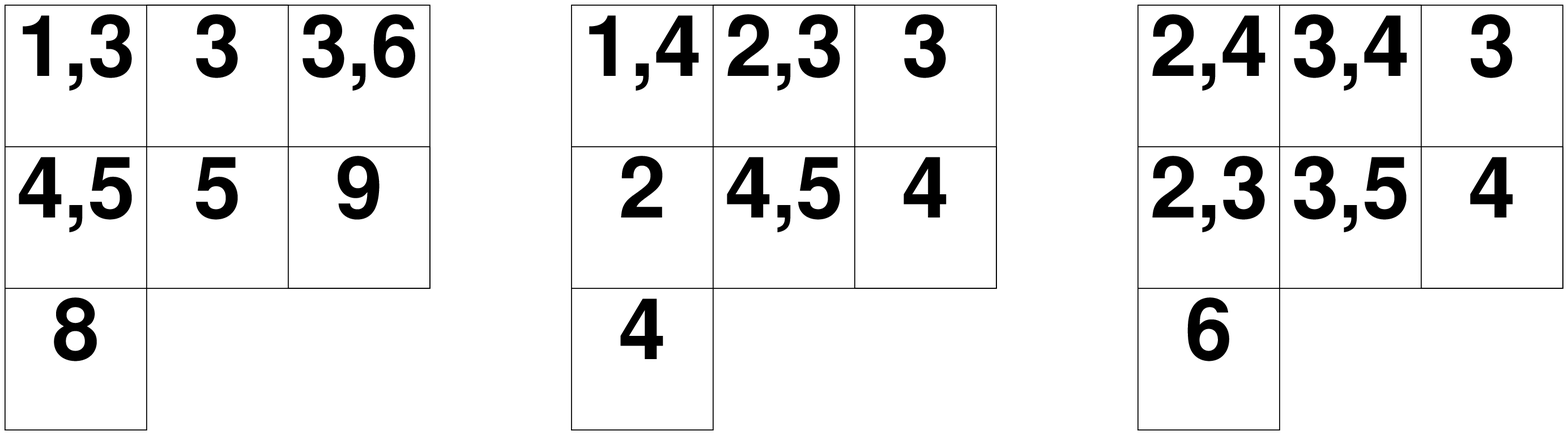,height=.7in}
\]
Hereafter, we will not bother to write the commas in our examples; no
confusion will result because we only use numbers that are at
most~$9$.

The \dfn{union} of two set-valued tableaux of the same shape $\lambda$
simply assigns to each box of $\lambda$ the union of the two sets
associated to it.  Moreover, if either set-valued tableau is limit
semistandard, so is the union.  The intersection is not always
defined, however, because of the requirement that every box
of~$\lambda$ be nonempty.

In addition to the partition $\lambda$, fix a maximum entry value $n
\in \naturals_+ = \{1,2,3,\ldots\}$.  Define the \dfn{empty-face
tableau} $E_{\lambda,n}$ associated to $\lambda$ and $n$ as the union
of all the semistandard tableaux with shape $\lambda$ and all entries
at most $n$.

\begin{figure}[tb]
  \centering
  \epsfig{file=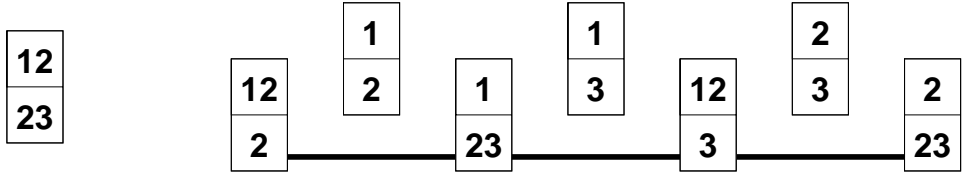,width=5in}
  \caption{A Young tableau complex.  At left is the empty-face tableau.}
  \label{fig:1ball}
\end{figure}

Consider the partition~$\lambda$ as a poset in which $(i,j) \leq
(i',j')$ whenever \mbox{$i \leq i'$} and \mbox{$j \leq j'$}; thus each
box is less than the boxes southeast of it.  Writing $[n] =
\{1,\ldots,n\}$, we get a tableau complex $\DELTA
\lambda{[n]}T{E_{\lambda,n}}$, in the sense of
Theorem~\ref{t:posetball}: take $\Psi$ to be the set of pairs
\mbox{$($upper box, lower box$)$} in which one box sits atop another
in~$\lambda$, so $T$ is the set of semistandard Young tableaux
on~$\lambda$ with maximum value~$n$.  Observing that $\DELTA
\lambda{[n]}T{E_{\lambda,n}}$ depends only on~$\lambda$ and~$n$, we
denote this \dfn{Young tableau complex} by $\Delta(\lambda,n)$.  See
Figure~\ref{fig:1ball} for an example.
The special case of Theorem~\ref{t:posetball} for Young tableau
complexes is as follows.

\begin{Corollary*}
  The Young tableau complex $\Delta(\lambda,n)$ is homeomorphic to a
  shellable ball or sphere, and its interior faces are labeled by
  Buch's semistandard set-valued Young tableaux \cite{buch:KLR}.
\end{Corollary*}

\begin{Example} \label{exa:21}
Let $\lambda = (2,1)$ and $n=3$.  Then $\Delta_{\lambda,3}$ is a
$3$-dimensional ball.  It has one interior vertex, missing the $2$ in
the upper left box.  We draw the boundary $2$-sphere in
Figure~\ref{fig:2sphere}.
\begin{figure}[htbp]
  \centering
  \epsfig{file=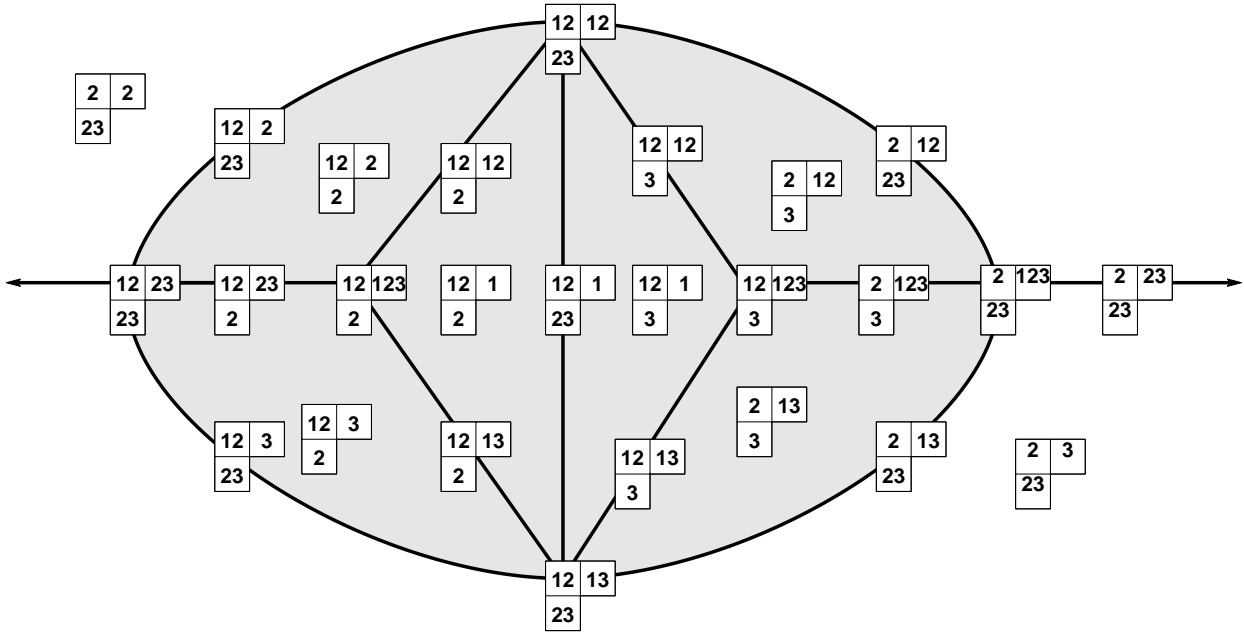,width=6.5in}
  \caption{A triangulated $2$-sphere of properly limit semistandard
  tableaux.  The two edges with arrows meet around the back side of
  the sphere.}
  \label{fig:2sphere}
\end{figure}%
\end{Example}%

\begin{figure}[htbp]
  \centering
  \epsfig{file=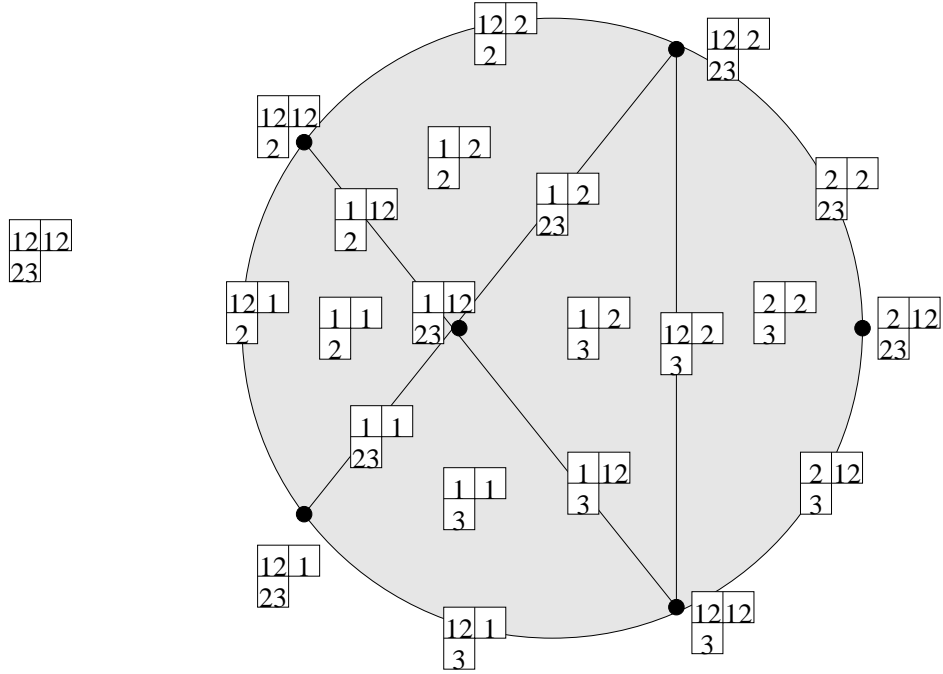,height=3.5in}
  \caption{The empty-face tableau and simplicial complex for $\lambda
    = (2,1)$ and $\vec n = (2,3)$.}
  \label{fig:23ball}
\end{figure}

When vertex-decomposing a Young tableau complex---that is, writing it
as the union of the star and the deletion of a single vertex---the two
subcomplexes are \dfn{flagged} Young tableau complexes, in which a
vector $\vec n$ bounds the sizes of the entries in the rows
of~$\lambda$.  This suggests that we ought to work in that level of
generality; see Figure~\ref{fig:23ball}.

Flagged Young tableaux are used to compute the Schubert polynomials
for vexillary permutations \cite{Wachs}, and one choice of vertex
decomposition parallels the ``transition formula'' for their double
Grothendieck polynomials~\cite{lascoux, lascoux:ASM}.  Hence we are
able to give set-valued-tableaux-based formulae for double
Grothendieck polynomials of vexillary permutations, via
Theorem~\ref{t:kpoly}.  The second formula in the following corollary,
which appeared already in \cite{Knutson.Miller.Yong}, is a common
generalization of Buch's and Wachs's formulae, each of which
specializes to the usual tableau formula for Schur polynomials.  The
other two parts give new formulae for these polynomials.  See
Section~\ref{sec:vex} for proofs.

\begin{Corollary*}
  Let $\pi\in S_n$ be a vexillary permutation with associated
  partition $\lambda$ and flagging $\vec n$.  Each of the following is
  a formula for the double Grothendieck polynomial $\groth_\pi({\bf
  x}, {\bf y})$.
  \begin{enumerate}
  \item
  As a sum over the set $\LSVT(\lambda,\vec n)$ of limit semistandard
  tableau associated to $(\lambda,\vec n)$,
  \[
    \groth_\pi({\bf x}, {\bf y}) = 
    \sum_{\tau \in \LSVT(\lambda,\vec n)} 
    \prod_{b \in \lambda}\prod_{i \in \tau(b)} (1-x_i y_{i+j(b)}^{-1}) 
    \prod_{h \in E_{\lambda,\vec n}(b) \setminus \tau(b)} 
    x_{h}y_{h+j(b)}^{-1},
  \]
  where $E_{\lambda,\vec n}=\BT$ is the union of all
  semistandard tableaux $\tau \in \SSYT(\lambda, \vec n)$,
  and $j(b) = c(b)-r(b)$ is the difference of the row and column
  indices of the box~$b\in\lambda$.

  \item
  As a sum over the set $\SVT(\lambda,\vec n)$ of semistandard
  set-valued tableau associated to $(\lambda,\vec n)$,
  \[
    \groth_\pi({\bf x}, {\bf y}) = 
    \sum_{\tau \in \SVT(\lambda,\vec n)} (-1)^{|\tau|-|\lambda|}
    \prod_{b \in \lambda} \prod_{i \in \tau(b)} (1-x_i y_{i+j(b)}^{-1}).
  \]
  
  \item
  As a sum over the set $\SSYT(\lambda,\vec n)$ of semistandard Young
  tableaux of shape $\lambda$ flagged by~$\vec n$,
  \[
    \groth_\pi({\bf x}, {\bf y}) = 
    \sum_{\tau \in \SSYT(\lambda, \vec n)} \prod_{b\in\lambda} \prod_{i\in
    \tau(b)} (1-x_{i} y_{i+j(b)}^{-1}) \prod_{h\in E_{\lambda, \vec n}(b)
    \setminus N_\tau(b)} x_{h}y_{h+j(b)}^{-1},
  \] 
  where $N_\tau$ is the tableau obtained by adding to each box $b$ all
  numbers in $E_{\lambda,\vec n}(b)$ either smaller than the entry of
  $\tau(b)$, or larger than that entry provided that replacing the entry
  with the larger number would not give a tableau in $\SSYT(\lambda,
  \vec n)$.
  \end{enumerate}
\end{Corollary*}

The second of these formulae for vexillary double Grothendieck
polynomials was based on the algebraic geometry of matrix Schubert
varieties \cite{Knutson.Miller.Yong}.  It was that geometry that first
motivated us to fit Young tableaux into a simplicial complex.

\section{Properties of tableau complexes}

\subsection{Generalities and boundary faces}\label{sub:gen}
Recall that the vertices of a tableau complex consist of the
complements $\bxty$ of single elements of~$E$.

\begin{Proposition}\label{prop:easy}
  Let $\Delta = \DELTA XYTE$ be a tableau complex, and assume that $(x
  \mapsto y) \in E$.
  \begin{enumerate}
  \item\label{easy:pure}
  $\Delta$ is pure, and its facets are labeled by the tableaux in~$T$.

  \item\label{easy:codim}
  Writing $|F| = \sum_{x \in X} |F(x)|$, the codimension of a face~$F$
  in $\DELTA XYTE$ equals $|F| - |X|$.

  \item\label{easy:shellable}
  Each ridge is contained in at most two facets.  In particular, if
  $\DELTA XYTE$ is shellable then it is homeomorphic to a ball or
  sphere.

  \item\label{easy:cone}
  $\bxty$ is a \dfn{cone vertex} (meaning it lies in every facet) if
  and only if $f(x)\neq y$ for every $f\in T$.  In particular, $E =
  \BT$ exactly when $\DELTA XYTE$ has no cone vertices.

  \item
  $\bxty$ is a \dfn{phantom vertex}, meaning $\bxty \notin \Delta$,
  precisely if $f(x)=y$ for all~\mbox{$f\in T$}.
  \end{enumerate}
\end{Proposition}

\begin{proof}
%
Only statement~\ref{easy:shellable} is not immediate from the
definitions.  A ridge is a set-valued tableau taking one extra value.
By the pigeonhole principle, since every $x \in X$ gets at least one
$y \in Y$, there exists exactly one~$x$ with two values from~$Y$, all
others being~$1$.  Such a set-valued tableau can contain at most two
tableaux from~$T$.  The statement about being a ball or sphere now
follows from~\cite[Proposition~4.7.22]{BLSWZ}.
%
\end{proof}

Cone vertices are in some sense uninteresting: a simplicial complex
can be canonically reconstructed from its set of cone vertices and its
\dfn{core}, which is the subcomplex with the cone vertices removed.
In particular, the whole complex is a ball or sphere if and only if
its core is a ball or sphere.
It is convenient for inductive purposes not to assume that $E = \BT$,
although we will generally assume $E = \BT$ in examples.

\begin{Proposition}\label{prop:boundary}
  Assume that $\DELTA XYTE$ is homeomorphic to a ball or sphere.  A
  face~$F$ of $\DELTA XYTE$ lies on the boundary of $\DELTA XYTE$ if
  and only if there exists a tableau $g: X \to Y$ such that $g
  \subseteq F$ but $g \notin T$.
\end{Proposition}
\begin{proof}
By definition, a face of a simplicial ball lies in the boundary if and
only if it is a face of a boundary ridge.  A ridge itself lies in the
boundary if and only if it is a face of precisely one facet.  Of the
two tableaux contained in a given ridge, at least one must lie in~$T$,
because a ridge is a face of $\DELTA XYTE$.  Hence a ridge is a
boundary face if and only if the unique other tableau it contains does
not lie in~$T$.

Now let $F$ be an arbitrary face of~$\DELTA XYTE$.  If every function
$f \subseteq F$ lies in~$T$, then every ridge with $F$ as a face is a
union of two tableaux from~$T$, so $F$ is interior.  On the other
hand, suppose that $g \subseteq F$ for some tableau $g \notin T$, and
let $f \in T$ be a facet having $F$ as a face, so $f \subseteq F$.
Then $f \cup g \subseteq F$ is a set-valued tableau that has $F$ as a
face.  Some of the elements $x \in X$ are assigned two distinct
$Y$-values by $f \cup g$.  If deleting the value $f(x)$ from $(f \cup
g)(x)$ yields a face~$G$ of $\DELTA XYTE$, then induction on the
codimension implies that $G$ lies on the boundary, and hence $F
\supseteq G$ does, as well.  If no such $x$ exists, so deleting the
value $f(x)$ from $f \cup g$ always results in a set-valued tableau
that is not a face of $\DELTA XYTE$, then $f$ is the unique tableau
in~$T$ with $f \cup g$ as a face; thus $f \cup g$ is a face of only
one facet (namely~$f$), and hence $f \cup g$ is a boundary face
with~$F$ as a subface.%
\end{proof}

\subsection{Safe vertices in tableau complexes}\label{sub:safe}

Given a simplicial complex $\Delta$ with a vertex~$v$, define the
\dfn{star} and \dfn{deletion} of~$v$ to be
\[
  \starr_v \Delta = \{C \in \Delta: C\cup {v} \in \Delta\} \quad
  \text{and}\quad \del_v \Delta = \{C\in \Delta : v\notin C\}.
\]
Then $\Delta = \starr_v \Delta \cup \del_v \Delta$.  The star has an
obvious cone vertex, namely $v$ itself, and its deletion from the star
is called the \dfn{link} of~$v$ in $\Delta$.  More generally, the
link of a face~$C$ in a simplicial complex~$\Delta$ is defined as
\[
  \link_C\Delta = \{D\in\Delta : D\cap C=\nothing, D\cup C\in\Delta\}.
\]
By convention, the vertex set of this link does not include the (now
phantom) vertices~of~$C$.

\begin{Proposition}\label{prop:link}
Let $F$ be a face of $\Delta = \DELTA XYTE$.  Let $T_\link = \{f \in
T: f \subseteq F\}$ be the set of facets of~$\Delta$ having $F$ as a
subface.  Then the link of~$F$ in~$\Delta$ is isomorphic to $\DELTA
XY{T_\link}F$.
\end{Proposition}
\begin{proof}
It follows from the definitions that the faces of both $\link_F\Delta$
and $\DELTA XY{T_\link}F$ are the set-valued tableaux contained in~$F$
and containing a tableau from~$T$.%
\end{proof}

\begin{Proposition}\label{prop:safevertex1}
Let $\DELTA XYTE$ be a tableau complex.  Let $T_\starr = \{f \in T:
f(x)\neq y\}$.  Then $\starr_\bxty \DELTA XYTE = \DELTA
XY{T_\starr}E$.
\end{Proposition}
\begin{proof}
Since $\DELTA XYTE$ is pure, the star of $\bxty$ is the union of the
(closures of) facets that have $\bxty$ as a vertex.  These facets are
exactly the tableaux $f\in T_\starr$.%
\end{proof}

Call a vertex $\bxty$ of~$\DELTA XYTE$ \dfn{safe} if for every $f\in
T$, changing the label on~$x$ from $f(x)$ to~$y$ yields a tableau that
is again in $T$.
While the star of a vertex in a pure complex is always pure, the
deletion might not be.
\begin{Proposition}\label{prop:safevertex2}
The deletion $\del_\bxty \Delta$ of the vertex $\bxty$ from the
simplicial complex~$\Delta = \DELTA XYTE$ is pure if and only if
either $\bxty$ is a cone vertex or $\bxty$ is safe.
\end{Proposition}
\begin{proof}
If $\bxty$ is a cone vertex then $\Delta$ is the cone over $\del_\bxty
\Delta$.  A simplicial complex is pure if and only if the cone over it
is, so we assume that $\bxty$ is not a cone vertex.

Given a set-valued tableau~$F$, let $\del_\bxty F$ denote the
set-valued tableau that sends $a \mapsto F(a)$ for $a\neq x$ and sends
$x \mapsto F(x) \cup \{y\}$.  In particular, $\del_\bxty F = F$ if and
only if $y \in F(x)$.  The definitions imply that $\del_\bxty \Delta$
consists of the set-valued tableaux $\del_\bxty F$ for $F \in \Delta$,
and the facets of $\del_\bxty F$ have the form $\del_\bxty f$ for
tableaux $f \in T$.  Since $\bxty$ is not a cone vertex, at least one
facet of~$\Delta$ sends $x$ to~$y$.  Thus the deletion is pure if and
only if, for all tableaux $f \in T$, the set-valued tableau
$\del_\bxty f$ contains a tableau $g \in T$ satisfying $g(x) = y$.
The desired result follows because when $f(x)$ does not already
equal~$y$, the only possibility for~$g$ is obtained
by changing
$f(x)$ to~$y$.%
\end{proof}

\begin{Corollary}\label{cor:safevertex}
Let $\bxty$ be a safe vertex of the tableau complex $\Delta = \DELTA
XYTE$.  If $T_\del = \{f \in T: f(x)=y\}$, then $\del_\bxty \Delta =
\DELTA XY{T_\del}E$.\qed
\end{Corollary}

\subsection{Tableau complexes on posets}

At this point we make some additional assumptions to guarantee a ready
supply of safe vertices.  The following theorem is stated much more
generally than our motivating examples require; we hope that this
Bourbakiesque level of generality helps to indicate which assumptions
are leading to which conclusions.

The key to our geometric conclusions (shellable ball or sphere) is the
notion of \dfn{vertex-decomposable} simplicial complex in the sense
of \cite{BP}.  By definition,
every simplex is vertex-decomposable, and an arbitrary simplicial
complex
is vertex-decomposable if and only if it is pure and has a vertex
whose deletion and link are both vertex-decomposable.

\begin{Lemma}[\cite{BP}]\label{l:shellable}
A simplicial complex~$\Delta$ is shellable if both the deletion
$\del_v\Delta$ and the star $\starr_v\Delta$ of a vertex $v$ are
shellable.  Hence all vertex-decomposable simplicial complexes are
shellable.
\end{Lemma}
\begin{proof}
(\cite{BP}) Construct a shelling of $\Delta$ by concatenating
shellings of the deletion and star of~$v$ (in that order), the latter
being the cone over a shelling of the link.%
\end{proof}

\begin{Theorem}\label{thm:posetball}
Let $X$ and $Y$ be finite partially ordered sets.  For each $x \in X$,
let $Y_x$ be a totally ordered subset of~$Y$.  Fix a set~$\Psi$ of
pairs $(x_1,x_2)$ from~$X$ such that $x_1 < x_2$.  Let $T$ be the set
of tableaux $f: X \to Y$ such that
\begin{itemize}
\item $f(x) \in Y_x$ for all $x\in X$;
\item $f$ is weakly order-preserving, 
  i.e.\ $x_1\leq x_2$ implies that $f(x_1)\leq f(x_2)$; and
\item if $(x_1,x_2)\in \Psi$ then $f(x_1) < f(x_2)$.
\end{itemize}
Let $E \subseteq X \times Y$ contain $\BT$.  Then $\DELTA XYTE$ is
homeomorphic to a ball or sphere, and it is vertex-decomposable.%
\end{Theorem}
\begin{proof}
We need only prove vertex-decomposability, for then the ball or sphere
conclusion is a consequence of
Proposition~\ref{prop:easy}.\ref{easy:shellable} and
Lemma~\ref{l:shellable}.  To demonstrate vertex-decomposa\-bility, we
need only find, for each tableau complex satisfying the hypotheses of
the theorem, a vertex whose deletion and link both satisfy the
hypotheses.

Suppose that $\DELTA XYTE$ has a cone vertex $\bxty$.  Viewing $\bxty$
as a subset of \mbox{$X \times Y$}, we find that $\bxty$ already
contains~$\BT$ by Proposition~\ref{prop:easy}.\ref{easy:cone}, so
$\DELTA XYT\bxty$ satisfies the conditions of the theorem.  This
simplicial complex is the link by Proposition~\ref{prop:link}, and it
equals the deletion because $\bxty$ is a cone vertex.

Now assume that $\DELTA XYTE$ has no cone vertices.  If all of the
vertices of~$\DELTA XYTE$ are phantom, then there is only one facet
and we are done.  Otherwise, there exists a non-phantom vertex
\mbox{$(m \notmapsto y)$}.  Choose one with maximal possible~$m$, and
let $y_m$ be the maximum element of~$Y_m$.  Since there are no cone
points, the values of all tableaux in~$T$ are fixed at elements $x >
m$: for each $x > m$ and all $f \in T$ there is some $y_x \in Y$ such
that $f(x) = y_x$.  Therefore, as $\mym$ is itself not a cone vertex,
we get $y_m \leq y_x$ for all $x > m$, and $y_m < y_x$ if $(m,x) \in
\Psi$.  It follows that the vertex $\mym$ is safe: we can safely
change the label on~$m$ from $f(m)$ to~$y_m$ to get another tableau
satisfying the three conditions to be in~$T$ because $y_m \geq f(m)$
for all $f \in T$.

(We used that $Y_m$ has a maximum element for this, but not that it is
totally ordered.)

Most of the work has now been done in Section~\ref{sub:safe}: if
$T_\starr = \{f\in T: f(m) \neq y_m\}$ and $T_\del = \{f\in T: f(m) =
y_m\}$, then the star and deletion of $\mym$ are $\DELTA
XY{T_\starr}E$ and $\DELTA XY{T_\del}E$, respectively, by
Proposition~\ref{prop:safevertex1} and Corollary~\ref{cor:safevertex}.
The star and deletion satisfy the three conditions from the statement
of the theorem, with the same $X$, $Y$, and~$\Psi$, but with $Y_m$
changed either to $Y_m \setminus \{y_m\}$ or else to $\{y_m\}$,
respectively.  Given that the star satisfies the hypotheses of the
theorem, arguing as in the second paragraph of the proof shows that
the link does, as well.

(To work inductively, we need $Y_m \setminus \{y_m\}$ to again have a
maximum element; this is why we required $Y_m$ to be totally ordered.
In addition, our new choice of~$m$ for the link must have a maximum
element in its~$Y_m$; this is why we need {\em every} $Y_x$ to
have a maximum.)%
\end{proof}

Before this theorem, we never needed to compare $f(x_1)$ and $f(x_2)$
for $x_1\neq x_2$; in some sense, it would have been more natural
for the tableaux to take values in separate sets $Y_x$. Now that
we used a partial order on $Y$ to define our set $T$ of tableaux,
we have finally made such comparisons.

\begin{Example}\label{exa:skew_shape}
More generally than in Section~1.2, let $X$ be the set of boxes in a
skew-shape $\lambda/\mu$, and each $Y_x = Y=\{1,\ldots,n\}$.
Partially order $X$ by asking that each box is less than the boxes
southeast of it.  Let $\Psi$ be the set of pairs $\{$(upper box, lower
box)$\}$ where one box is atop another.  Then $T$ is the set of
semistandard Young tableaux of shape $\lambda/\mu$ with maximum value
$n$, and $\DELTA XYT{\bigcup_{f \in T} f}$ is the \dfn{Young skew
tableau complex}.

The faces of this complex are labeled with {\em set-valued Young skew
tableaux}, which were also introduced in \cite{buch:KLR}.  Buch's
definition of ``semistandard'' set-valued Young tableaux exactly
matches our criterion, Proposition~\ref{prop:boundary}, for a face to
be interior.

(In fact the $\Psi$ machinery was unnecessary to model
semistandardness; we could just take $\Psi=\nothing$, subtract $r-1$
from the values in the $r$th row, and adjust the sets $Y_x$ to get a
set combinatorially equivalent to semistandard Young tableaux.  But
the formulation with $\Psi$ is clearer, more general, and no more
difficult.)
\end{Example}

\begin{Example}
  Let $X$ be a poset, $\Psi=\nothing$, and $Y=\{0,1\}$. Then the
  tableaux correspond to partitioning $X$ into a lower and an upper
  order ideal (the $0$ and $1$ parts), or equivalently to antichains
  in $X$ (the maximum elements labeled $0$). By
  Theorem~\ref{thm:posetball}, this tableau complex is homeomorphic to
  a ball (or sphere, if $X$ is totally unordered).  
\end{Example}

\begin{Remark} \label{rk:greedoid}
Other classes of vertex-decomposable complexes include the greedoid
complexes \cite{BKL} and subword complexes \cite{KM:subword}; see
\cite[Remark~2.6]{KM:subword} for an extended discussion.  Tableau
complexes are different from each of these.  For example, the Young
tableau complex for the vertical domino with entries at most~$5$ is
not a greedoid complex if the ground set is taken to be the vertex
set.  To show the difference between subword and tableau complexes,
consider the Young tableau complex for the $2 \times 2$ square shape
with entries at most~$3$; it has dimension~$3$ and eight vertices,
none of which is a cone vertex.
On the other hand, deleting all cone points from the subword complex
in \cite{Knutson.Miller.Yong} having the same tableaux for facets
yields a simplicial complex of dimension~$2$ with seven~vertices.
\end{Remark}

It is worth noting that the phrase ``ball or sphere'' essentially
always really means ``ball''.  To get a sphere, there must be no cone
vertices, so $E = \BT$.  But even then, every ridge lies in two
facets, so every vertex must be safe; in other words, the possible
$T$-tableau values at each $x \in X$ are independent.  We spell this
out further in Section~\ref{sec:extremal}.  For now, here is a
characterization of the interior faces, which includes all of the
faces in the case of a sphere.  As a matter of notation, if $Y_1$ and
$Y_2$ are two subsets of a poset~$Y$, write $Y_1 \leq Y_2$ if $y_1
\leq y_2$ for all $y_1 \in Y_1$ and $y_2 \in Y_2$.  Similarly, write
$Y_1 < Y_2$ if strict inequality holds.  The following is an immediate
consequence of the definitions and Proposition~\ref{prop:boundary}.

\begin{Corollary}\label{c:boundary}
Assume the notation and hypotheses of Theorem~\ref{thm:posetball}.  A
face $F$ is interior to $\DELTA XYTE$ if and only if
  \begin{itemize}
  \item
  $F(x_1) \leq F(x_2)$ whenever $x_1 < x_2$; and
  \item
  $F(x_1) < F(x_2)$ whenever $(x_1 < x_2) \in \Psi$.
\end{itemize}
\end{Corollary}

\subsection{Shelling poset tableau complexes}

The next theorem will help us describe the $h$-vectors and Hilbert
series of poset tableau complexes and their Stanley-Reisner rings.

\begin{Theorem}\label{thm:shelling}
Assume the notation and hypotheses of Theorem~\ref{thm:posetball}, and
choose a linear extension $\varepsilon$ of the partial ordering
on~$X$.  Lexicographically order $T$ by comparing $f_1,f_2 \in T$ at
the $\varepsilon$-largest element $m \in X$ where they differ.
Placing the one with the larger label on~$m$ first yields a total
order on the facets of $\DELTA XYTE$ that is a shelling.
\end{Theorem}
\begin{proof}
Remember that $f_1(m)$ and $f_2(m)$ are comparable, since $Y_m$ is
totally ordered.  Therefore the procedure in the statement of the
theorem yields a total order on the facets.  We will show that this
total order is the shelling produced by applying
Lemma~\ref{l:shellable} recursively as in the proof of
Theorem~\ref{thm:posetball}.  At each stage in that proof, we either
vertex decompose at a cone vertex or we choose a maximal element $m
\in X$ among those supporting non-phantom vertices $(m \notmapsto y)$.
Vertex decomposing at a cone vertex does not alter the set~$T$ of
facet tableaux, so it does not matter in which order we delete cone
vertices.  Only the order in which we choose the maximal elements~$m$
matters.  Use the $\varepsilon$-order: since Lemma~\ref{l:shellable}
puts the deletion ($f(m) = y_m$) first before the star ($f(m) \neq
y_m$, which is equivalent to $f(m) < y_m$ because $y_m$ is maximum
in~$Y_m$), the resulting shelling is as desired.%
\end{proof}

\section{Characterizations of tableau complexes}\label{sec:extremal}

Let $\Delta$ be a pure simplicial complex on a vertex set~$V$.
Declare that $W \subseteq V$ is a \dfn{pure factor} of~$\Delta$ if the
number $|f \cap W|$ of vertices in the intersection of~$f$ with~$W$ is
the same for all facets $f \in \Delta$.  For example, a singleton
$\{v\}$ is a pure factor if and only if $v$ is a phantom or cone
vertex, with $|f \cap \{v\}| = 0$ or $|f \cap \{v\}| = 1$,
respectively.  If $W$ is a pure factor, then its complement $V
\setminus W$ is a pure factor, too.  For any set $W \subseteq V$ of
vertices, write $\Delta|_W$ for the \dfn{full subcomplex} $\del_{V
\setminus W} \Delta$ \dfn{supported on~$W$}.

\begin{Proposition}\label{prop:join}
Let $\Delta$ be a pure simplicial complex on the vertex set $V$, and
suppose that $V = V_1 \cup \cdots \cup V_k$ is partitioned into a
disjoint union of pure factors $V_1, \ldots, V_k$.  Then
$\Delta|_{V_i}$ is pure for each $1 \leq i \leq k$, and $\Delta$ is a
top-dimensional subcomplex of their \dfn{join}
\[
  \{F_1 \cup \cdots \cup F_k : F_i \text{ is a face of }
  \Delta|_{V_i}\, \text{ for each } 1 \leq i\leq k\}.
\]
\end{Proposition}
\begin{proof}
Since $\Delta$ is pure, it follows from the definitions that
$\Delta|_{V_i}$ is pure.  On the other hand, it also follows by
definition that a subset $f \subseteq V_1 \cup \cdots \cup V_k$ is a
facet of the join of the complexes $\Delta|_{V_i}$ if and only if $f
\cap V_i$ is a facet of the individual complex $\Delta|_{V_i}$ for
each~$i$.  Since the vertex set~$V$ is the disjoint union $V_1 \cup
\cdots \cup V_k$, every facet of~$\Delta$ has this property.%
\end{proof}

\begin{Lemma}\label{l:purefactor}
  Let $\Delta$ be a tableau complex $\DELTA XYT\BT$.  For each $x\in
  X$, define the subset $V_x = \{\bxty: (x\mapsto y) \in \BT\}$ of the
  vertex set of~$\Delta$.  Then
\begin{enumerate}
\item
  the subsets $\{V_x : x \in X\}$ partition the vertex set of~$\Delta$;

\item\label{I}
  each subset $V_x$ is a pure factor; in fact, $|f \cap V_x| =
  |V_x|-1$ for every facet tableau~$f$; and

\item\label{II}
  each induced complex $\Delta|_{V_x}$ is the boundary of the
  simplex on~$V_x$.
\end{enumerate}
\end{Lemma}
\begin{proof}
The first two numbered claims are immediate from the definitions.  For
the third,
it follows from the second that $\Delta|_{V_x}$ is a union of some
subset of the facets (each of size $|V_x|-1$) in the boundary of the
simplex on~$V_x$.  Each facet of $\Delta|_{V_x}$ avoids using some
(unique) vertex of~$V_x$.  If any vertex of $V_x$ does not occur this
way, then it lies in every facet of $\Delta|_{V_x}$; in other words,
it is a cone vertex of~$\Delta|_{V_x}$.  Since each facet of~$\Delta$
survives after deleting $V \setminus V_x$ to give a facet of
$\Delta|_{V_x}$, we conclude that $\Delta$ has a cone vertex,
contradicting Proposition~\ref{prop:easy}.\ref{easy:cone}.  Therefore
every face of size $|V_x| - 1$ occurs in $\Delta|_{V_x}$.%
\end{proof}
  
\begin{Theorem} \label{thm:pure}
A pure complex is (isomorphic to) a tableau complex $\DELTA XYT\BT$ if
and only if it can be expressed as a top-dimensional subcomplex of a
join of boundaries of simplices.
\end{Theorem}
\begin{proof}
That $\DELTA XYT\BT$ can be expressed in the desired manner is an
immediate consequence of Proposition~\ref{prop:join} and
Lemma~\ref{l:purefactor}.  Now suppose that that $\Delta$ is a pure
complex expressible as a top-dimensional subcomplex of the join of
boundaries of simplices with vertex sets $V_1,\ldots,V_k$.  Then the
vertex set $V$ of~$\Delta$ is the disjoint union $V_1 \cup \cdots \cup
V_k$.  Let $X = \{1,\ldots,k\}$ and set $Y = V$.  Each facet $f$
of~$\Delta$ defines a function $X \to Y$ taking $i$ to the element
$V_i \setminus f$.  Using these as the set $T$ of tableaux, we find
that $\Delta = \DELTA XYT\BT$.%
\end{proof}

\begin{Remark} \label{rk:sphere}
In particular, if a tableau complex has no boundary ridges (ridges
contained in just one facet), then it is the join of a bunch of
boundaries of simplices, and in particular it is a sphere.  This, plus
Proposition~\ref{prop:link}, gives another proof of
Proposition~\ref{prop:boundary}.
\end{Remark}

Let the \dfn{codimension} of a pure complex $\Delta$ be the number of
vertices outside any facet.  The only way for the codimension to
equal~$0$ is if $V$ consists only of cone vertices (so $V$ is a
simplex).  If $V$ breaks up as a union $V_1\cup\cdots\cup V_k$
of pure factors, then the codimension
of $\Delta$ is the sum of the codimensions of the full subcomplexes
supported on the $V_i$.

This suggests a characterization of tableau complexes by the following
extremal property.  Given a pure complex $\Delta$ with the cone and
phantom vertices deleted, look for a pure factor~$W$.  Splitting into
$W$ and $V\setminus W$, the codimension of each full subcomplex can be
no larger than that of~$\Delta$.  By the theorem, $\Delta$ is a
tableau complex if and only if we can split enough to whittle the
codimensions of all of the full subcomplexes down to~$1$.

The situation is somewhat dual to order complexes of ranked posets.
If~$P$ is a ranked poset, its \dfn{order complex} has vertex set~$P$,
and $Q\subset P$ defines a face if and only if $Q$ is totally ordered.
If~$P_r$ denotes the set of elements with a given rank~$r$, then the
induced complex on~$P_r$ is pure of dimension~$0$, rather than
codimension~$1$ like a tableau complex.
(If it seems unsatisfying for ``codimension~$1$'' to be dual to
``dimension~$0$'', then consider the latter more honestly as
``affine-dimension~$1$''.)  Very few simplicial complexes are both
order complexes and tableau complexes; we leave their characterization
as an exercise for the reader.

\begin{Remark} \label{rk:matroid}
Tableau complexes bear superficial similarities to matroid complexes.
A simplicial complex is a matroid if and only if every subcomplex
induced on a subset of the vertex set is pure.  Theorem~\ref{thm:pure}
says that a simplicial complex is a tableau complex if and only if the
vertex set can be partitioned into subsets that are pure factors of
codimension~1.  In reality, there are matroid complexes that are not
tableau complexes, and conversely.  For example, we have already seen
in Remark~\ref{rk:greedoid} that tableau complexes can fail to be
greedoid complexes, of which matroid complexes are special cases.  For
an example the other way, the matroid for the complete graph $K_4$ on
four vertices is the union of three segments joined at a point.  If
this were a tableau complex, then so would be the result of deleting
the cone point, by Proposition~\ref{prop:link}.  But a set of three
points is not a tableau complex by Remark~\ref{rk:sphere}, and hence
neither is the matroid for~$K_4$.
\end{Remark}

\section{Hilbert series and \K-polynomials}\label{sec:hilb}

In this section we collect some formulae for the Hilbert series of the
Stanley-Reisner rings of tableau complexes.  Our main reference for
generalities on Betti numbers, Hilbert series, and \K-polynomials
(which are numerators of Hilbert series) is \cite{cca}.  For notation,
let $\kk$ be a field, and set $S = \kk[V]$, the polynomial ring in
variables $v \in V$ indexed by a finite set~$V$.  This is the ambient
ring for objects like the \dfn{Stanley-Reisner ideal} $\ID =
\<\prod_{v \in F} v : F$ is not a face of~$\Delta\>$ of a simplicial
complex~$\Delta$ with vertex set~$V$, and the \dfn{Stanley-Reisner
ring} $S/\ID$.  We shall use the alphabet $\tt = \{t_v : v \in
V\}$ for finely graded Hilbert series and \K-polynomials.  When
$\Delta$ is a tableau complex $\DELTA XYTE$, recall that $V$ is the
set $\{\bxty : (x \mapsto y) \in E\}$ of complements of single
elements of~$E$.

\begin{Proposition} \label{p:kpoly1}
The \K-polynomial associated to the tableau complex $\Delta = \DELTA
XYTE$ is
\[
  K(S/\ID;\tt) = \sum_F \prod_{x \in X} \Bigg(\prod_{y \in E(x)
  \setminus F(x)} t_\bxty \prod_{y \in F(x)} (1-t_\bxty)\Bigg),
\]
the sum being over all set-valued tableaux $F \subseteq E$ each
containing some tableau $f \in T$.
\end{Proposition}
\begin{proof}
This formula is \cite[Theorem~1.13]{cca} applied to tableau complexes,
since the condition $y \in E(x) \setminus F(x)$ means that $\bxty$ is
a vertex of~$F$ and the condition $y \in F(x)$ means that $\bxty$ is
not a vertex of~$F$.%
\end{proof}

Our second formula uses the ball-or-sphere hypothesis; it therefore
holds for (among other things) all poset tableau complexes.  It will
be simpler to prove the formula for general balls and spheres first.

\begin{Proposition}\label{p:Hilb2}
If $\Delta$ is a simplicial ball or sphere with vertex set~$V$, then
\[
  K(S/\ID;\tt) = \sum_F (-1)^{\codim(F)}\prod_{v \in V \setminus
  F}(1-t_v),
\]
where the sum is over all interior faces of~$\Delta$.  (All faces are
interior if $\Delta$ is a sphere.)%
\end{Proposition}
\begin{proof}
Consider the \dfn{Alexander dual} ideal $\IDstar = \<\prod_{v \in V
\setminus F} v : F$ is a face of~$\Delta\>$, and start by calculating
the \K-polynomial of~$\IDstar$.  The coefficient on the monomial
$\prod_{v \in V \setminus F} t_v$ in $K(\IDstar;\tt)$ is the
alternating sum of the Betti numbers of~$\IDstar$ in degree $\prod_{v
\in V \setminus F} v$ \cite[Proposition~8.23]{cca}.  By Hochster's
formula \cite[Corollary~1.40]{cca}, the $i^{\rm th}$ such Betti number
equals the dimension $\dim_\kk \widetilde H_{i-1}(\link_F \Delta;\kk)$
of the reduced homology of the link of~$F$ in~$\Delta$, and it comes
with a sign $(-1)^i$.  If $F$ is a boundary face, then the link of~$F$
is contractible; but if $F$ is an interior face, then the link of~$F$
is a sphere of dimension $\codim(F) - 1$.  Therefore $K(\IDstar;\tt) =
\sum_F (-1)^{\codim(F)} \prod_{v \in V \setminus F} t_v$, where the
sum is over all interior faces~$F$ of~$\Delta$.  The Alexander
inversion formula \cite[Theorem~5.14]{cca} now implies the desired
result.%
\end{proof}

\begin{Theorem} \label{t:Kball}
If the tableau complex $\Delta = \DELTA XYTE$ is homeomorphic to a
ball or sphere, then
\[
  K(S/\ID;\tt) = \sum_F (-1)^{|F|-|X|} \prod_{x \in X} \prod_{y \in
  F(x)} (1-t_\bxty),
\]
the sum being over all set-valued tableaux $F \subseteq E$ such that
every tableau $f \subseteq F$ lies in~$T$.
\end{Theorem}
\begin{proof}
The factor $(-1)^{|F|-|X|}$ is the codimension of~$F$.  The condition
$y \in F(x)$ means that $\bxty$ lies in the vertex set of~$\Delta$ but
not~$F$.  The sum is over all interior faces by
Proposition~\ref{prop:boundary}.  Therefore the result is simply
Proposition~\ref{p:Hilb2} for tableau complexes.%
\end{proof}

A \dfn{shelling} of a pure $d$-dimensional simplicial complex~$\Delta$
is an ordering of the facets $F_1,\ldots, F_k$ such that $F_i \cap
(F_1 \cup \cdots \cup F_{i-1})$ has pure dimension $d-1$ for each $1
\leq i \leq k$.  This guarantees that for each~$i$, there is a unique
minimal new face $N_i \in \Delta$ that is a face of~$F_i$ but not of
$F_1,\ldots,F_{i-1}$.  By convention, $N_1$ is the empty face.

\begin{Lemma}\label{l:Ni}
Given a shelling of a simplicial complex~$\Delta$ with new faces
$N_1,\ldots,N_k$ as above,
\[
  K(S/\ID;\tt) = \sum_{i=1}^{k}\prod_{v\not\in F_i}(1-t_v)\prod_{v\in
  N_i} t_v.
\]
\end{Lemma}
\begin{proof}
Use induction on the number~$k$ of facets of~$\Delta$.%
\end{proof}

The ${\mathbb Z}$-graded coarsening of the Hilbert series to one
variable~$t$ gives
\[
  H(S/\ID,t) = \sum_{j=0}^{d}\frac{h_j t^j}{(1-t)^d}.
\]
When $\Delta$ is shellable, the \dfn{$h$-vector}
$(h_0,h_1,\ldots,h_d)$ consists of nonnegative integers.  Moreover,
the shelling gives a manifestly positive way to compute these numbers:
$h_j$ counts the number of dimension~$j$ faces among $N_1,\ldots,
N_k$.

\begin{Theorem}\label{thm:kpoly}
Resume the notation and hypotheses of Theorem~\ref{thm:posetball}, and
set $E = \BT$.  Given a tableau $f \in T$, define $U_f$ as the set of
elements $y \in Y$ such that $f(x) \leq y$ and moving the label on~$x$
from $f(x)$ up to~$y$ yields a tableau in~$T$.  Then
\[
  K(S/\ID;\tt) = \sum_{f\in T} \prod_{x\in X}\bigg(\big(1-t_{(x
  \notmapsto f(x))}\big) \prod_{y \in U_f(x)} t_\bxty\bigg).
\]
Finally, if $\eta(f) = -|X| + \sum_{x\in X} |U_f(x)|$, then $h_j$ is
the number of tableaux $f \in T$ with $\eta(f) = j$.
\end{Theorem}
\begin{proof}
The proof will be done once we produce a shelling of~$\Delta$ for
which $N_f = E(x) \setminus U_f(x)$ is the minimal new face at the
stage when we add the facet~$f$.  Pick a linear extension
$\varepsilon$ of~$X$ and take the shelling order of the facets $f_1,
f_2, \ldots$ of~$\Delta$ as in Theorem~\ref{thm:shelling}.  For
$f:=f_i$ we show that $N_f$ is the minimal new face of~$f$.
  
First, $N_f$ is a set-valued tableau in~$\Delta$ that is a face
of~$f$, since it contains~$f$.  Second, to see that $N_f$ is not a
face of any previous facet, we must show that $f$ does not contain
$f_1,\ldots, f_{i-1}$.  Note that by construction, any other $g \in T$
contained in $N_f$ must assign to each $x \in X$ either $f(x)$ or some
$y < f(x)$.  Such a facet tableau $g$ must appear later than~$f$ in
the facet ordering.  The maximality of $N_f \subseteq E$
containing~$f$ and not containing $f_1,\ldots,f_{i-1}$ is also clear
from the construction.  Hence $N_f$ is the minimal new face of~$f$, as
claimed.  For the \K-polynomial formula, apply Lemma~\ref{l:Ni}.%
\end{proof}

\begin{Example}
For the tableau complex in Figure~\ref{fig:23ball} (after
Example~\ref{exa:21}), listing the facets in the order 
\[
\begin{array}{c}
  \tableau{2&2\\
           3}\,,\,\,
  \tableau{1&2\\
           3}\,,\,\,
  \tableau{1&1\\
           3}\,,\,\,
  \tableau{1&2\\
           2}\,,\,\,
  \tableau{1&1\\
           2}
  \\[2ex]
\end{array}
\]
yields the shelling in the proof of Theorem~\ref{thm:kpoly}.  For the
first of these tableaux, all of the sets $U(x)$ are singletons: there
is no way to increase the number in any box while respecting the
flagging, which requires that the entries in the top row are at
most~$2$ and the lower entry is at most~$3$.  For the second tableau
above, $U(x) = \{1,2\}$ for the upper-left corner~$x$, because moving
the~$1$ up to a~$2$ keeps the tableau semistandard.  Similarly, all of
the sets $U(x)$ for the third and fourth tableaux are singletons
except for upper-right box and the bottom box, respectively, whose
sets $U(x)$ are $\{1,2\}$ and $\{2,3\}$.  In the last tableau above,
only the two lower-right corners have non-singleton sets $U(x)$, and
these are $\{1,2\}$ and $\{2,3\}$.  For the above five tableau, the
function $\eta$ at the end of Theorem~\ref{thm:kpoly} takes the values
$0$, $1$, $1$, $1$, and $2$, respectively.  Thus the simplicial
complex in Figure~\ref{fig:23ball} has $h$-vector $(1,3,1)$.
\end{Example}

Our final \K-polynomial formula in this section will arise again after
Corollary~\ref{cor:formulas}.

\begin{Proposition} \label{p:kpoly2}
If $\bxty$ is a safe vertex of $\Delta = \DELTA XYTE$, then
\[
  K(S/\ID;\tt) = t_\bxty K(S/I_\del;\tt) + (1-t_\bxty)
  K(S/I_\starr;\tt),
\]
where $I_\del$ and $I_\starr$ are the Stanley-Reisner ideals for the
deletion tableau complex $\DELTA{X}Y{T_\del}E$ and the star tableau
complex $\DELTA{X}Y{T_\starr}E$ from
Propositions~\ref{prop:safevertex1} and~\ref{prop:safevertex2},
respectively.
\end{Proposition}
\begin{proof}
Any vertex decomposition $\Delta = \del_v \Delta \cup \starr_v \Delta$
gives an inductive formula
\[
  K(S/I_\Delta;\tt) = t_v K(S/I_{\del_v \Delta};\tt) + (1 - t_v)
  K(S/I_{\starr_v \Delta};\tt)
\]
for the \K-polynomial.%
\end{proof}

\section{Applications to vexillary double Grothendieck polynomials}\label{sec:vex}

In this section, we apply Section~\ref{sec:hilb} to obtain formulae
for double Grothendieck polynomials for vexillary permutations.  This
gives formulae for the Hilbert series and \K-polynomials of vexillary
matrix Schubert varieties (also known as (one-sided) ladder
determinantal varieties).  See \cite{Knutson.Miller.Yong} for a
treatment of the related algebraic geometry.

\subsection{Vexillary permutations and flaggings of partitions}
Identify a permutation \mbox{$\pi \in S_n$} with the square array
having blank boxes in all locations except at $(i,\pi(i))$ for $i =
1,\ldots,n$, where we place \dfn{dots}.  This defines the
\dfn{dot-matrix} of $\pi$.  We associate the \dfn{diagram}
\[
  D(\pi)=\big\{(p,q) \in \{1,\ldots,n\} \times \{1,\ldots,n\} :
  \pi(p)>q \hbox{ and } \pi^{-1}(q)>p\big\}
\]
to~$\pi$.  Pictorially, if we draw a ``hook'' consisting of lines
going east and south from each dot, then $D(\pi)$ consists of the
squares not in the hook of any dot.

\begin{Example}\label{example:diag}
Let $\pi=\left(\begin{array}{ccccccccc}
1 & 2 & 3 & 4 & 5 & 6 & 7 & 8 & 9 \\
8 & 7 & 1 & 6 & 2 & 9 & 5 & 3 & 4
\end{array}\right)$.   Its dot-matrix and diagram are combined below:
\[
\begin{picture}(120,120)
\put(0,0){\makebox[0pt][l]{\framebox(120,120)}}
\put(100,113.33){\circle*{4}}
\thicklines
\put(100,113.33){\line(1,0){20}}
\put(100,113.33){\line(0,-1){113.33}}
\put(86.67,100){\circle*{4}}
\put(86.67,100){\line(1,0){33.33}}
\put(86.67,100){\line(0,-1){100}}
\put(6.67,86.67){\circle*{4}}
\put(6.67,86.67){\line(1,0){113.33}}
\put(6.67,86.67){\line(0,-1){86.67}}

\put(73.33,73.33){\circle*{4}}
\put(73.33,73.33){\line(1,0){46.67}}
\put(73.33,73.33){\line(0,-1){73.33}}

\put(20,60){\circle*{4}}
\put(20,60){\line(1,0){100}}
\put(20,60){\line(0,-1){60}}

\put(113.33,46.67){\circle*{4}}
\put(113.33,46.67){\line(1,0){6.67}}
\put(113.33,46.67){\line(0,-1){46.67}}

\put(60,33.33){\circle*{4}}
\put(60,33.33){\line(1,0){60}}
\put(60,33.33){\line(0,-1){33.33}}
\put(33.33,20){\circle*{4}}
\put(33.33,20){\line(1,0){86.67}}
\put(33.33,20){\line(0,-1){20}}

\put(46.67,6.67){\circle*{4}}
\put(46.67,6.67){\line(1,0){73.33}}
\put(46.67,6.67){\line(0,-1){6.67}}

\thinlines
\put(0,106.67){\makebox[0pt][l]{\framebox(93.33,13.33)}}
\put(0,93.33){\makebox[0pt][l]{\framebox(80,26.67)}}
\put(13.33,93.33){\line(0,1){26.67}}
\put(26.67,93.33){\line(0,1){26.67}}
\put(40,93.33){\line(0,1){26.67}}
\put(53.33,93.33){\line(0,1){26.67}}
\put(66.67,93.33){\line(0,1){26.67}}

\put(13.33,66.67){\makebox[0pt][l]{\framebox(53.33,13.33)}}
\put(26.67,66.67){\line(0,1){13.33}}
\put(40,66.67){\line(0,1){13.33}}
\put(53.33,66.67){\line(0,1){13.33}}

\put(26.67,26.67){\makebox[0pt][l]{\framebox(26.67,26.67)}}
\put(26.67,40){\line(1,0){40}}
\put(66.67,40){\line(0,1){13.5}}
\put(53.43,53.58){\line(1,0){13.34}}
\put(40,26.67){\line(0,1){26.67}}
\end{picture}
\]
\end{Example}

In what follows, we will assume our permutations $\pi$ are
\dfn{vexillary}, also known as \dfn{$2143$-avoiding}: there
exist no indices $1\leq a<b<c<d\leq n$ with
\mbox{$\pi(b)<\pi(a)<\pi(d)<\pi(c)$}.  We need some facts about
vexillary permutations; further details consistent with the
terminology and notation used here may be found
in~\cite{Knutson.Miller.Yong} and the references therein.

Throughout we will identify a partition $\lambda$ with its Young
diagram.  There is a partition $\lambda$ associated to $\pi$: let the
$k$th diagonal of $\lambda$ (those boxes $\{ (i,k+i)\}$) have as many
boxes as the $k$th diagonal of $D(\pi)$, for each $k$.  Indeed, this
sets up a natural bijection between the boxes of $\lambda$ and the
boxes of $D(\pi)$, taking the $j$th box down in the $k$th diagonal to
the $j$th box down in the $k$th diagonal.  (The difference is that in
$D(\pi)$ there may be spaces in between the boxes.)  This bijection
also defines a \dfn{flagging} $\vec n$ on the rows of $\lambda$.
Namely, $n_i\in {\mathbb N}_{+}$ equals the row of $D(\pi)$ containing
the box corresponding to the rightmost box of the $i^{th}$ row of
$\lambda$.  We will thus speak interchangeably about $\pi$ and the
pair $(\lambda,\vec n)$.

In Example~\ref{example:diag}, the permutation $\pi$ is vexillary,
$\lambda=(7,6,4,3,2)$, and $\vec n = (1,2,4,6,7)$.

We remark that $\vec n$ need not be a weakly increasing sequence.  For
instance,
if $\sigma=\left(\begin{array}{cccccccc}
1 & 2 & 3 & 4 & 5 & 6 & 7 & 8 \\
2 & 7 & 4 & 5 & 8 & 1 & 3 & 6 
\end{array}\right)$,
then $\lambda=(5,3,2,2,1)$ and $\vec n = (2,5,4,4,5)$.

Call a set-valued tableau $\tau$ with shape $\lambda$ \dfn{$\vec
n$-flagged} if the maximum (so, the last) entry in each row is bounded
above by the corresponding entry of $\vec n$.

Extend the definition of the \dfn{empty-face tableau} in the obvious
way: it is the union of all the $\vec n$-flagged semistandard tableau
on the shape $\lambda$.  Let this set-valued tableau be denoted by
$E_{\lambda, {\vec n}}$.  (Note that $E_{\lambda,\vec n}(b)$ is an
interval in the natural numbers~$\naturals$: the smallest entry is the
row position of $b\in \lambda$ while the largest entry is the position
of the corresponding box of $D(\pi)$, under the bijection between
$\lambda$ and $D(\pi)$ described above.)  This gives rise to a tableau
complex $\Delta(\lambda, \vec n)$ generalizing that described in
Section~1.2:

\begin{Corollary}\label{cor:flagged}
For a partition $\lambda$ and a flagging $\vec n$ associated to a
vexillary permutation~$\pi$, $\Delta(\lambda,\vec n)$ is a simplicial
ball, and its interior faces are the flagged semistandard set-valued
Young~tableaux.
\end{Corollary}

\subsection{Formulae for vexillary Grothendieck polynomials}
For each permutation $\pi\in S_n$ there is a \dfn{(double)
Grothendieck polynomial}
\[
  \groth_{\pi}(x_1,\ldots, x_n, y_1,\ldots y_n) \in \integers[x_1^{\pm
  1},\ldots,x_n^{\pm 1},y_1^{\pm 1},\ldots,y_n^{\pm 1}]
\]
of Lascoux and Sch\"utzenberger \cite{LS}.  The case that $\pi$ is
vexillary has been of specific interest; see \cite{fulton:duke,
KM:annals, Knutson.Miller.Yong}.  We give tableau formulae in this
setting.

Let $\SVT(\lambda,\vec n)$ denote the semistandard set-valued tableaux
of shape~$\lambda$ and flagging~$\vec n$.  Similarly, denote by
$\SSYT(\lambda, \vec n)$ and $\LSVT(\lambda,\vec n)$ the corresponding
set of semistandard and limit semistandard tableaux, respectively.
For a set-valued tableau~$\tau$, let $\tau(b)$ denote the set of
entries in box~$b$.

\begin{Corollary}\label{cor:formulas}
Let $\pi\in S_n$ be a vexillary permutation and $(\lambda,\vec n)$ be
the associated partition and flagging.  Each of the following is a
formula for the double Grothendieck polynomial $\groth_{\pi}({\bf x},
{\bf y})$.
\[
  \sum_{\tau\in \LSVT(\lambda,\vec n)} 
  \prod_{b\in \lambda}\prod_{i\in \tau(b)} (1-x_i y_{i+j(b)}^{-1}) 
  \prod_{i\in E_{\lambda,\vec n}(b)\setminus \tau(b)} 
  x_{i}y_{i+j(b)}^{-1}
\]
\[
  \sum_{\tau\in \SVT(\lambda,\vec n)} (-1)^{|\tau|-|\lambda|}
  \prod_{b\in\lambda} \prod_{i\in \tau(b)} (1-x_i y_{i+j(b)}^{-1})
\]
\[
  \sum_{\tau\in \SSYT(\lambda, \vec n)} \prod_{b\in\lambda}\prod_{i\in
  \tau(b)} (1-x_{i} y_{i+j(b)}^{-1})\prod_{i\in E_{\lambda, \vec
  n}(b)\setminus N_\tau(b)} x_{i}y_{i+j(b)}^{-1}
\]
Here, $j(b)=c(b)-r(b)$ is the difference of the column and row indices
of the box~$b\in\lambda$.  Moreover, in the last formula, $N_\tau$ is
the tableau obtained by adding to each box~$b$ all entries of~$E(b)$
either smaller than the entry of~$\tau(b)$, or larger than the entry
of~$\tau(b)$ but such that replacing $\tau(b)$ with this larger number
would not give a tableau in $\SVT(\lambda, \vec n)$.
\end{Corollary}
\begin{proof}
Formally, the second formula in the corollary is obtained as follows.
Compute the \K-polynomial of $\Delta_{\lambda,\vec n}$ via
Theorem~\ref{t:Kball}, using $\lambda$, $b$, and~$i$ here in place of
$X$, $x$, and~$y$ there.  Then, for each fixed $b$ and~$i$, substitute
the expression~$x_iy_{i+j(b)}^{-1}$ for the variable~$t_{(b \notmapsto
i)}$.

It was shown in~\cite{Knutson.Miller.Yong} that the second formula
equals the desired double Grothen\-dieck polynomial.  Therefore,
applying the same substitution procedure to the results of
Proposition~\ref{p:kpoly1} and Theorem~\ref{thm:kpoly} yields two more
formulae for the same Grothendieck polynomial.  These formulae are,
respectively, the first and third formulae here.%
\end{proof}

In the above proof we appealed to \cite{Knutson.Miller.Yong} to
confirm that our \K-polynomials are in fact Grothendieck polynomials.
Let us briefly sketch how this can be proved directly; we refer the
reader to \cite{Knutson.Miller.Yong} for terminology.  To each
vexillary permutation $\pi$ there is an \emph{accessible box} of
$\lambda$.  From this one can define two vexillary permutations
$\pi_{\mathcal P}$ and $\pi_{\mathcal C}$.  We obtain a safe vertex of
the flagged Young tableau complex by removing the largest entry of
$E_{\lambda, \vec n}$ that appears in the accessible box.  The
resulting deletion and star subcomplexes are naturally isomorphic to
(multicones over) the flagged Young tableau complexes for
$\pi_{\mathcal P}$ and $\pi_{\mathcal C}$ respectively.  The recursion
from Proposition~\ref{p:kpoly2} is precisely Lascoux's transition
formula for vexillary Grothendieck polynomials~\cite{lascoux,
lascoux:ASM} (after the substitution $t_\bxty \mapsto
x_{i}y_{i+j(b)}^{-1}$).  Thus, since both polynomials satisfy the same
recursion (and initial conditions), they are equal.

Specializations of these formulae are of interest.  Suppose we set
$y_j=1$ for each $j$ and replace $x_i$ with $1-x_i$ for each~$i$.  If
we assume furthermore that $\pi$ is Grassmannian, then we obtain
Buch's formula~\cite{buch:KLR} for the single Grothendieck polynomial
$\groth_{\lambda}({\bf 1-x})$ \cite{buch:KLR}.  If instead we take the
lowest degree terms of the polynomial, we obtain Wachs's formula for a
flagged Schur polynomial \cite{Wachs}.  Making both of these
specializations gives the classical tableau formula for an ordinary
Schur polynomial.

We remark that it is also possible to use similar methods to extend these
results to give set-valued skew tableau formulae for ``$321$-avoiding 
permutations'' (see Example~\ref{exa:skew_shape}).

\begin{excise}{%
    \comment{The crystal operations are discontinuous.  What else sucks?}
    
    \comment{A subword complex labeled by limit semistandard tableaux}
    
    \subsection{A ball with faces labeled by ordinary Young tableaux}
    
    For completeness, we include another construction of a simplicial
    ball with faces labeled by ordinary (not set-valued)
    tableaux.  This is based on a different Gr\"obner degeneration, of
    a Schubert variety in the Grassmannian in its Pl\"ucker embedding
    \cite{DEP}.
    
    Fix a partition $\lambda$ with $n$ boxes.  Define a \dfn{Hodge
    tableau} as a tableau $\tau$ on $\lambda$, with labels from
    $\{1,\ldots,n\}$ with the following restriction: if we use $k$,
    then the number of labels $j$ with $j<k$ must be exactly
    $k-1$.  Also, the labels should be {\em weakly} increasing {\em
    both} down and to the right.   One mnemonic: consider the tableau
    as the result of a race, where the boxes have run to the corner;
    the numbers are the ranks, where multiple ties are possible.
    
    If we let $\tau_k$ be the set of boxes with $\tau$-label at most
    $k$, then $\tau$ can be uniquely reconstructed from the weakly
    increasing sequence $(\tau_1,\tau_2,\ldots,\tau_n = \lambda)$.
    Conversely, from any strictly increasing sequence of partitions,
    with last entry $\lambda$, we can construct a corresponding
    $\tau$.
    
    This gives a rule for labeling each chain in the lattice of
    partitions contained in $\lambda$, and hence each simplex in the
    order complex of this lattice, by a Hodge tableau.  Such order
    complexes (and much more general ones) were proven in \cite{BW} to
    be homeomorphic to balls.  \comment{better check [DEP] to see if
    they already said this} In terms of the box-racing metaphor, one
    simplex is a face of another if the second tableau arises by
    breaking ties in the race of the first tableau.
    
    \begin{figure}[htbp]
      \centering
      \epsfig{file=Hodge.eps,width=4.5in}
      \caption{The simplicial ball for the $(3,2)$ partition, 
        minus its two cone vertices.}
      \label{fig:Hodge}
    \end{figure}
    
    Each vertex is labeled by a partition with only two different
    numbers ($1$ and something else), and a face is labeled by the
    maximum of the labels on all its vertices.  The facets are labeled
    by the standard tableaux with each number $\{1,\ldots,n\}$ used
    exactly once.
}\end{excise}%

\section*{Acknowledgements}
We would like to thank Eric Babson, Nantel Bergeron, Christian Krattenthaler, 
Carsten Lange and Ravi Vakil for helpful comments.
AK was supported by NSF grant DMS-0303523. 
EM was partially supported by NSF grant DMS-0304789, NSF CAREER grant
DMS-0449102 and a University of Minnesota McKnight Land-Grant
Professorship.
This work was partially completed while AY was an NSERC visitor to the
Fields Institute during the 2005 semester on ``The geometry of string
theory'', and while an NSF supported visitor at the Mittag-Leffler
institute, during the 2005 semester on ``Algebraic combinatorics''.
All of us would like to thank the 2005 AMS Summer Research
Institute on Algebraic Geometry in Seattle, 
where this work was also carried out.

\end{document}